\documentclass[12pt]{article}
\usepackage{latexsym}
 \usepackage{mathrsfs}
\usepackage{amssymb}
\usepackage{graphicx}

\newtheorem{Theorem}{Theorem}[part]
\newtheorem{Definition}{Definition}[part]
\newtheorem{Proposition}{Proposition}[part]

\newtheorem{Remark}{Remark}[part]
\newtheorem{Example}{Example}[part]
\newtheorem{Algorithm}{Algorithm}[part]

\def \ep{\hbox{ }\hfill$\Box$}

\def \d{\Delta}

\newcounter{RomanNumber}

\begin{document}
\title{Formulating an $n$-person noncooperative game as a tensor complementarity problem
}

\author{Zheng-Hai Huang\thanks{Department of Mathematics, School of Science, Tianjin University, Tianjin 300072, P.R. China. 
Email: huangzhenghai@tju.edu.cn. Tel: +86-22-27403615, Fax: +86-22-27403615. This author's work was supported by the National Natural Science Foundation of China (Grant No. 11431002).}
\and
Liqun Qi\thanks{Department of Applied Mathematics, The Hong Kong Polytechnic University, Hung Hom, Kowloon, Hong Kong, P.R. China. Email: liqun.qi@polyu.edu.hk. This author's work was supported by the Hong Kong Research Grant Council (Grant No. PolyU 502111, 501212, 501913 and 15302114).}
}

\date{February 10, 2016}

\maketitle

\begin{abstract}
\noindent

In this paper, we consider a class of $n$-person noncooperative games, where the utility function of every player is given by a homogeneous polynomial defined by the payoff tensor of that player, which is a natural extension of the bimatrix game where the utility function of every player is given by a quadratic form defined by the payoff matrix of that player. We will call such a problem the multilinear game. We reformulate the multilinear game as a tensor complementarity problem, a generalization of the linear complementarity problem; and show that finding a Nash equilibrium point of the multilinear game is equivalent to finding a solution of the resulted tensor complementarity problem. Especially, we present an explicit relationship between the solutions of the multilinear game and the tensor complementarity problem, which builds a bridge between these two classes of problems. We also apply a smoothing-type algorithm to solve the resulted tensor complementarity problem and give some preliminary numerical results for solving the multilinear games.
\vspace{3mm}

\noindent {\bf Key words:}\hspace{2mm} Game theory, $n$-person noncooperative game, Nash equilibrium, bimatrix game, tensor complementarity problem  \vspace{3mm}

\noindent {\bf Mathematics Subject Classifications (2000):}\hspace{2mm} 90C33, 65K10, 74P10, 91A10, 91A06. \vspace{3mm}

\end{abstract}

\section{Introduction}

\hspace{4mm} The $n$-person noncooperative game plays a fundamental yet important role in the development of game theory \cite{bo-95,or-94}. Nash \cite{nash-50,nash-51} proposed a very important concept of equilibrium, called Nash equilibrium, for $n$-person noncooperative games, which is a stable outcome in the sense that a unilateral deviation from a Nash equilibrium point by one of the players does not increase the payoff of that player. Nash \cite{nash-50,nash-51} has shown that every game of this kind has at least one equilibrium point in mixed strategies. The $n$-person noncooperative game and its various extensions have been extensively studied, for example, see \cite{bo-95,bomze-86,fk-07,gw-03,hyf-05,hf-report,kj-09,ll-13,Myerson-99,weibull-96,Yuan-11} and references therein.

A large number of economic models are formulated in terms of some $n$-person noncooperative game \cite{ah-92,fp-97}. In these applications, one of main concerns is how to find effectively a Nash equilibrium point, which depends heavily on the good mathematical description of the problem. For the two-person noncooperative game, one of the most popular models is the bimatrix game \cite{gs-96,lh-64}, where the utility function of every player is a quadratic form defined by the payoff matrix of that player. It is well known that the bimatrix game can be reformulated as a linear complementarity problem \cite{cps-92,lh-64,ms-64}. The first approach for finding Nash equilibrium point of a two-person game was proposed in \cite{lh-64}, which was designed based on the reformulated linear complementarity problem of the concerned game.

The polymatrix game is an important subclass of $n$-person noncooperative games, where the payoff of one player relative to the decisions of any other player is independent of the remaining players' choices \cite{Howson-72}. The utility function of each player is the sum of $n-1$ quadratic forms where every quadratic form is defined by the payoff matrix of this player with respect to any other player. Obviously, the polymatrix game is an extension of the bimatrix game. It is well known that the polymatrix game can be reformulated as a linear complementarity problem \cite{hxq-06,Howson-72}.

Recently, Song and Qi \cite{sq-15} introduced a class of complementarity problems, called tensor complementarity problems, where the involved function is defined by some homogeneous polynomial of degree $n$ with $n>2$. It is known that the tensor complementarity problem is a generalization of the linear complementarity problem \cite{cps-92}; and a subclass of nonlinear complementarity problems \cite{fp-03,hxq-06}. The tensor complementarity problem was studied recently by many scholars \cite{bhw-15-r,cqw-16,dlq-2015-r,glqx-15-r,hsw-15-r,lqx-15-r,sq-15-r,sq-15-r1,sq-15-r2,sy-15-r,whb-15-r}.

In this paper, we consider a class of $n$-person noncooperative games, where the utility function of every player is a homogeneous polynomial of degree $n$ defined by the payoff tensor of that player. The new model is a natural extension of the bimatrix game; and we call it the multilinear game in this paper. We will reformulate the multilinear game as a tensor complementarity problem. We show that finding a Nash equilibrium point of the multilinear game is equivalent to finding a solution of the resulted tensor complementarity problem; and especially, we exhibit an explicit corresponding relation between the solutions of these two classes of problems. In addition, we also apply a smoothing-type algorithm to solve the resulted tensor complementarity problem and give some preliminary numerical results for solving some multilinear games.

\section{Preliminaries}
\hspace{4mm} In this section, we introduce some notation and give some basic results, which will be used in the subsequent analysis.

Throughout this paper, we assume that $m_1, m_2, \cdots, m_n$ and $n$ are positive integers, and $n>2$ unless specifically stated.
For any positive integer $n$, we denote $\{1,2,\ldots,n\}$ by $[n]$ and the $n$-dimensional vector of ones by $e_n$.

A real $n$-th order $m_1\times m_2\times \cdots \times m_n$-dimensional tensor $\mathscr{B}$ is a multiple array in $\mathbb{R}^{m_1\times m_2\times\cdots \times m_n}$, which can be written as $\mathscr{B}:=(b_{i_1i_2\cdots i_n})$ where $b_{i_1i_2\cdots i_n}\in \mathbb{R}$ for any $i_j\in [m_j]$ and $j\in [n]$. If $m_1=m_2=\cdots=m_n=l$, then $\mathscr{B}$ is called a real $n$-th order $l$-dimensional tensor. We will denote the set of all real $n$-th order $l$-dimensional tensors by $\mathbb{T}_{n,l}$.

We will use the following concept, which can be found in \cite{KB-09}.
\begin{Definition}\label{def-product}
The $k$-mode (vector) product of a tensor $\mathscr{B}=(b_{i_1i_2\cdots i_n})\in \mathbb{R}^{m_1\times m_2\times\cdots \times m_n}$ with a vector $v\in \mathbb{R}^{m_k}$ is denoted by $\mathscr{B}\bar{\times}_k v$, which is a real $(n-1)$-th order $m_1\times\cdots m_{k-1}\times m_{k+1}\times \cdots \times m_n$-dimensional tensor with
$$
(\mathscr{B}\bar{\times}_k v)_{i_1\cdots i_{k-1}i_{k+1}\cdots i_n}=\sum_{i_k=1}^{m_k}b_{i_1i_2\cdots i_n}v_{i_k}
$$
for any $i_j\in [m_j]$ with $j\in [n]\setminus \{k\}$.
\end{Definition}

For any tensor $\mathscr{B}=(b_{i_1i_2\cdots i_n})\in \mathbb{R}^{m_1\times m_2\times\cdots \times m_n}$ and any vector $u^k\in \mathbb{R}^{m_k}$ with $k\in [n]$, we will use $\mathscr{B}u^1u^2\cdots u^n$ to denote $\mathscr{B}\bar{\times}_1 u^1\bar{\times}_2 u^2\bar{\times}_3\cdots \bar{\times}_n u^n$ and use $\mathscr{B}u^2\cdots u^n$ to denote $\mathscr{B}\bar{\times}_2 u^2\bar{\times}_3\cdots \bar{\times}_n u^n$ for simplicity. Then, by using Definition \ref{def-product}, we have
$$
\mathscr{B}u^1u^2\cdots u^n=\sum\limits_{i_1=1}^{m_1}\sum\limits_{i_2=1}^{m_2}\cdots \sum\limits_{i_n=1}^{m_n}b_{i_1i_2\cdots i_n} u^1_{i_1}u^2_{i_2}\cdots u^n_{i_n}
$$
and
$$
\mathscr{B}u^{2}\cdots u^{n}=\left(\begin{array}{c}
\sum\limits_{i_2=1}^{m_2}\cdots\sum\limits_{i_n=1}^{m_n} b_{1i_2\cdots i_n} u^2_{i_2}\cdots u^n_{i_n}\\ \vdots \\ \sum\limits_{i_2=1}^{m_2}\cdots\sum\limits_{i_n=1}^{m_n} b_{m_1i_2\cdots i_n} u^2_{i_2}\cdots u^n_{i_n}
\end{array}\right).
$$

For any $k\in [n]$, it is easy to see that
\begin{eqnarray*}
& &\frac{\partial}{\partial u^k}\mathscr{B}u^1u^2\cdots u^n\\
& &\quad =\left(\begin{array}{c}
\sum\limits_{i_1=1}^{m_1}\cdots\sum\limits_{i_{k-1}=1}^{m_{k-1}}\sum\limits_{i_{k+1}=1}^{m_{k+1}}\cdots\sum\limits_{i_n=1}^{m_n} b_{i_1\cdots i_{k-1}1i_{k+1}\cdots i_n} u^1_{i_1}\cdots u^{k-1}_{i_{k-1}}u^{k+1}_{i_{k+1}}\cdots u^n_{i_n}\\ \vdots \\ \sum\limits_{i_1=1}^{m_1}\cdots\sum\limits_{i_{k-1}=1}^{m_{k-1}}\sum\limits_{i_{k+1}=1}^{m_{k+1}}\cdots\sum\limits_{i_n=1}^{m_n} b_{i_1\cdots i_{k-1}m_ki_{k+1}\cdots i_n} u^1_{i_1}\cdots u^{k-1}_{i_{k-1}}u^{k+1}_{i_{k+1}}\cdots u^n_{i_n}
\end{array}\right)\\
& &\quad =\left(\begin{array}{c}
\sum\limits_{i_1=1}^{m_1}\cdots\sum\limits_{i_{k-1}=1}^{m_{k-1}}\sum\limits_{i_{k+1}=1}^{m_{k+1}}\cdots\sum\limits_{i_n=1}^{m_n} \bar{b}_{1i_1\cdots i_{k-1}i_{k+1}\cdots i_n} u^1_{i_1}\cdots u^{k-1}_{i_{k-1}}u^{k+1}_{i_{k+1}}\cdots u^n_{i_n}\\ \vdots \\ \sum\limits_{i_1=1}^{m_1}\cdots\sum\limits_{i_{k-1}=1}^{m_{k-1}}\sum\limits_{i_{k+1}=1}^{m_{k+1}}\cdots\sum\limits_{i_n=1}^{m_n} \bar{b}_{m_ki_1\cdots i_{k-1}i_{k+1}\cdots i_n} u^1_{i_1}\cdots u^{k-1}_{i_{k-1}}u^{k+1}_{i_{k+1}}\cdots u^n_{i_n}
\end{array}\right).
\end{eqnarray*}

We introduce the following tensors.
\begin{Definition}\label{def-B-bar}
For any tensor $\mathscr{B}=(b_{i_1i_2\cdots i_n})\in \mathbb{R}^{m_1\times m_2\times\cdots \times m_n}$ and any $k\in [n]$, we define tensor
$$
\bar{\mathscr{B}}^k:=(\bar{b}_{i_1i_2\cdots i_n})\in \mathbb{R}^{m_k\times m_1\times \cdots \times m_{k-1}\times m_{k+1}\times \cdots\times m_n}
$$
with
$$
\bar{b}_{i_1i_2\cdots i_n}=b_{i_ki_1\cdots i_{k-1}i_{k+1}\cdots i_n},\quad \forall i_j\in [m_j]\;\; \mbox{\rm and}\;\; j\in [n].
$$
\end{Definition}

Then, the following results can be easily obtained.
\begin{Proposition}\label{prop-1}
For any tensor $\mathscr{B}=(b_{i_1i_2\cdots i_n})\in \mathbb{R}^{m_1\times m_2\times\cdots \times m_n}$ and any $k\in [n]$, let $\bar{\mathscr{B}}^k$ be defined by Definition \ref{def-B-bar}. Then, $\bar{\mathscr{B}}^1=\mathscr{B}$ and any $k\in [n]$,
$$
\frac{\partial}{\partial u^k}\mathscr{B}u^1u^2\cdots u^n=\bar{\mathscr{B}}^ku^1\cdots u^{k-1}u^{k+1}\cdots u^n
$$
and
$$
\langle u^k,\bar{\mathscr{B}}^ku^1\cdots u^{k-1}u^{k+1}\cdots u^n\rangle=\mathscr{B}u^1u^2\cdots u^n
$$
hold for any $k\in [n]$.
\end{Proposition}

For any tensor $\mathscr{B}=(b_{i_1i_2\cdots i_n})\in \mathbb{T}_{n,m}$ and any vector $u\in \mathbb{R}^m$, we will use $\mathscr{B}u^{n-1}$ to denote the vector $\mathscr{B}\bar{\times}_2 u\bar{\times}_3\cdots \bar{\times}_n u\in \mathbb{R}^m$ for simplicity. Then, by using Definition \ref{def-product}, we have
\begin{eqnarray}\label{E-bum1}
\mathscr{B}u^{n-1}=\left(\begin{array}{c}
\sum\limits_{i_2=1}^{m}\cdots\sum\limits_{i_n=1}^{m} b_{1i_2\cdots i_n} u_{i_2}\cdots u_{i_n}\\ \vdots \\ \sum\limits_{i_2=1}^{m}\cdots\sum\limits_{i_n=1}^{m} b_{mi_2\cdots i_n} u_{i_2}\cdots u_{i_n}
\end{array}\right).
\end{eqnarray}
In fact, such a notation has been extensively used in the literature \cite{Qi-05}.

In the following, we denote $m:=\sum_{j=1}^nm_j$. We will use $x=\left((x^k)_{k\in [n]}\right)\in \mathbb{R}^m$ and $x^*=\left(({x^k}^*)_{k\in [n]}\right)\in \mathbb{R}^m$ to mean that
$$
x=\left(\begin{array}{c} x^1 \\ x^2 \\ \vdots \\ x^n\end{array}\right),\; x^*=\left(\begin{array}{c} {x^1}^* \\ {x^2}^* \\ \vdots \\ {x^n}^*\end{array}\right)\in \mathbb{R}^{m_1} \times \mathbb{R}^{m_2}\times\cdots \times \mathbb{R}^{m_n}=\mathbb{R}^m.
$$

\section{Description of the multilinear game}
\hspace{4mm}
The so-called multilinear game is a noncooperative game with a finite number of players, each with a finite number of pure strategies, which is specified as follows.
\begin{itemize}
\item[(I)] There are $n$ players: player $1$, player $2$, $\cdots$, player $n$, i.e., the set of players is $[n]$.
\item[(II)] For any $k\in [n]$, player $k$ has $m_k$ pure strategies, i.e., the pure strategy set of player $k$ is $[m_k]$.
\item[(III)] For any $k\in [n]$, let $\mathscr{A}^k=(a^k_{i_1i_2\cdots i_n})$ be payoff tensor of player $k$, that is to say, for any $i_j\in [m_j]$ with any $j\in [n]$, if player $1$ plays his $i_1$-th pure strategy, player $2$ plays his $i_2$-th strategy, $\cdots$, player $n$ plays his $i_n$-th strategy, then the payoffs of player $1$, player $2$, $\cdots$, player $n$ are $a^1_{i_1i_2\cdots i_n}$, $a^2_{i_1i_2\cdots i_n}$, $\cdots$, $a^n_{i_1i_2\cdots i_n}$, respectively.
\item[(IV)] For any $k\in [n]$, let $x^k=(x_{i_j}^k)\in \mathbb{R}^{m_k}$ represent a mixed strategy of player $k$, where $x_{i_j}^k\geq0$ is the relative probability that player $k$ plays his $i_j$-th pure strategy for any $i_j\in [m_k]$, i.e., $x^k\in \Omega_k:=\{x\in \mathbb{R}^{m_k}: x\geq 0\; \mbox{\rm and}\; e_{m_k}^Tx = 1\}$.
\end{itemize}
Thus, the utility function of player $k$ is
\begin{eqnarray}\label{utility-func}
 \mathscr{A}^kx^1x^2\cdots x^n=\sum_{i_1=1}^{m_1}\sum_{i_2=1}^{m_2}\cdots\sum_{i_n=1}^{m_n} a^k_{i_1i_2\cdots i_n} x^1_{i_1}x^2_{i_2}\cdots x^n_{i_n}
\end{eqnarray}
for any $k\in [n]$.

We say that $x=\left((x^k)_{k\in [n]}\right)\in \mathbb{R}^m$ is a joint mixed strategy if $x^k$ is a mixed strategy of player $k$ for any $k\in [n]$, i.e., $x^k$ satisfies $x^k\geq 0$ and $e_{m_k}^Tx^k = 1$.

A joint mixed strategy $x^*=\left(({x^k}^*)_{k\in [n]}\right)\in \mathbb{R}^m$ is said to be a Nash equilibrium point of the multilinear game, if for any joint mixed strategy $x=\left((x^k)_{k\in [n]}\right)\in \mathbb{R}^m$ and any $k\in [n]$, it holds that
$$
\mathscr{A}^k{x^1}^*{x^2}^*\cdots {x^n}^*\geq \mathscr{A}^k{x^1}^*\cdots{x^{k-1}}^*x^k{x^{k+1}}^*\cdots {x^n}^*.
$$

It is obvious that a joint mixed strategy $x^*=\left(({x^k}^*)_{k\in [n]}\right)\in \mathbb{R}^m$ is a Nash equilibrium point of the multilinear game if and only if, for any $k\in [n]$, $x^*=\left(({x^k}^*)_{k\in [n]}\right)\in \mathbb{R}^m$ is an optimal solution of the following optimization problem:
\begin{eqnarray}\label{E-opt-1}
\begin{array}{cl}
  \max\limits_{x^k\in\mathbb{R}^{m_k}} & \mathscr{A}^k{x^1}^*\cdots{x^{k-1}}^*x^k{x^{k+1}}^*\cdots {x^n}^* \vspace{2mm} \\
  {\rm s.t.}  & e_{m_k}^Tx^k = 1, x^k \geq 0.
\end{array}
\end{eqnarray}

\begin{Remark}
\begin{itemize}
\item[(i)] The model (\ref{E-opt-1}) of the multilinear game is given by using the Nash equilibrium in mixed strategies, not in pure strategies. Throughout this paper, we consider such a model.
\item[(ii)] There are many different models on the $n$-person noncooperative game. One of general models is that the utility function is a continuously differentiable concave function and each set $\Omega_k\; (\forall k\in [n])$ defined in (IV) is convex.
\item[(iii)] If in (III), we let $A^{kj}$ denote the payoff matrix of player $k$ with respect to player $j$ (i.e., if player $k$ plays his $p$-th pure strategy and player $j$ plays his $q$-th pure strategy, then the payoff of player $k$ is $a_{pq}^{kj}$); and furthermore, instead of (\ref{utility-func}), we define the utility function of player $k$ by
    $$
    {x^k}^T\sum_{j\in [n]\setminus \{k\}}A^{kj}x^j,
    $$
    then the corresponding problem is the polymatrix game.
\item[(iv)] It is obvious that the multilinear game considered in this paper is different from the polymatrix game. Both the multilinear game and the polymatrix game are the generalizations of the bimatrix game, however, it seems that the multilinear game is a more natural extension of the bimatrix game than the polymatrix game.
\end{itemize}
\end{Remark}

Without loss of generality, we assume in this paper that $a^k_{i_1i_2\cdots i_n}>0$ for any $k\in [n]$ and any $i_j\in [m_j]$ with all $j\in [n]$. In fact, it is obvious that there exists a sufficient large $c>0$ such that $a^k_{i_1i_2\cdots i_n}+c>0$ for any $k\in [n]$ and any $i_j\in [m_j]$ with all $j\in [n]$. Since for any joint mixed strategy $x=\left((x^k)_{k\in [n]}\right)\in \mathbb{R}^m$ and any $k\in [n]$, we have that
$$
\sum_{i_1=1}^{m_1}\sum_{i_2=1}^{m_2}\cdots\sum_{i_n=1}^{m_n} (a^k_{i_1i_2\cdots i_n}+c) x^1_{i_1}x^2_{i_2}\cdots x^n_{i_n}=\mathscr{A}^kx^1x^2\cdots x^n+c,
$$
it is easy to see that $x^*=\left(({x^k}^*)_{k\in [n]}\right)\in \mathbb{R}^m$ is a Nash equilibrium point of the multilinear game with payoff tensors $\mathscr{A}^k$ for all $k\in [n]$ if and only if $x^*=\left(({x^k}^*)_{k\in [n]}\right)\in \mathbb{R}^m$ is a Nash equilibrium point of the multilinear game with payoff tensors $\mathscr{A}^k+c\mathscr{E}$  for all $k\in [n]$, where $\mathscr{E}\in {\mathbb R}^{m_1\times m_2\times \cdots \times m_n}$ is a tensor whose all entries are 1.

\section{Reformulation of the multilinear game}
\hspace{4mm} For any given tensor $\mathscr{B}\in \mathbb{T}_{n,l}$ and vector $q\in \mathbb{R}^l$, the tensor complementarity problem, denoted by the TCP$(q,\mathscr{B})$, is to find a vector $z\in \mathbb{R}^l$ such that
$$
z\geq 0,\quad \mathscr{B}z^{n-1}+q\geq 0, \quad \langle z, \mathscr{B}z^{n-1}+q\rangle=0,
$$
which was introduced recently by Song and Qi \cite{sq-15}; and was further studied by many scholars \cite{bhw-15-r,cqw-16,dlq-2015-r,glqx-15-r,hsw-15-r,lqx-15-r,sq-15-r,sq-15-r1,sq-15-r2,sy-15-r,whb-15-r}. When $n=2$, the tensor $\mathscr{B}$ reduces to a matrix, denoted by $B$; and the TCP$(q,\mathscr{B})$ becomes: find a vector $z\in \mathbb{R}^l$ such that
$$
z\geq 0,\quad Bz+q\geq 0, \quad \langle z, Bz+q\rangle=0,
$$
which is just the linear complementarity problem \cite{cps-92}.

In this section, we show that the multilinear game can be reformulated as a specific tensor complementarity problem.

Using payoff tensors $\mathscr{A}^k$ for all $k\in [n]$, we construct a new tensor:
\begin{eqnarray}\label{E-new-tensor}
\mathscr{A}:=(a_{i_1i_2\cdots i_n})\in \mathbb{T}_{n,m}
\end{eqnarray}
where
\begin{eqnarray*}
a_{i_1i_2\cdots i_n}=\left\{\begin{array}{l}
a^1_{i_1(i_2-m_1)\cdots (i_n-\sum_{j=1}^{n-1}m_j)},\\ \qquad \mbox{\rm if}\;\; i_1\in [m_1], i_2\in [m_1+m_2]\setminus [m_1], \cdots, i_n\in [\sum_{j=1}^nm_j]\setminus [\sum_{j=1}^{n-1}m_j], \vspace{2mm}\\
a^2_{(i_1-m_1)i_2(i_3-m_1-m_2)\cdots (i_n-\sum_{j=1}^{n-1}m_j)},\\ \qquad \mbox{\rm if}\;\; i_1\in [m_1+m_2]\setminus [m_1], i_2\in [m_1], \\
\qquad\quad i_3\in [\sum_{j=1}^3m_j]\setminus [m_1+m_2],\cdots, i_n\in [\sum_{j=1}^nm_j]\setminus [\sum_{j=1}^{n-1}m_j], \vspace{2mm}\\
a^k_{(i_1-\sum_{j=1}^{k-1}m_j)(i_2-m_1)\cdots (i_{k-1}-\sum_{j=1}^{k-2}m_j)i_k(i_{k+1}-\sum_{j+1}^km_j)\cdots (i_n-\sum_{j=1}^{n-1}m_j)}, \\ \qquad \mbox{\rm if}\;\; k\in [n]\setminus\{1,2\},\; \mbox{\rm and for any given}\; k, i_1\in [\sum_{j=1}^km_j]\setminus [\sum_{j=1}^{k-1}m_j], \\ \qquad\quad i_2\in [m_1+m_2]\setminus [m_1], \cdots, i_{k-1}\in [\sum_{j=1}^{k-1}m_j]\setminus [\sum_{j=1}^{k-2}m_j], i_k\in [m_1],\\ \qquad\quad i_{k+1}\in [\sum_{j=1}^{k+1}m_j]\setminus [\sum_{j=1}^{k}m_j], \cdots, i_n\in [\sum_{j=1}^nm_j]\setminus [\sum_{j=1}^{n-1}m_j],\vspace{2mm}\\
0, \quad \mbox{\rm otherwise}
\end{array}\right.
\end{eqnarray*}
for any $i_j\in [m]$ with $j\in [n]$.

For convenience of description, we introduce the following tensors by using the payoff tensors.
\begin{Definition}\label{def-tensor-bar-a}
For any $k\in [n]$, let $\mathscr{A}^k$ be the payoff tensor of player $k$; and define
$$
\bar{\mathscr{A}}^k:=(\bar{a}^k_{i_1i_2\cdots i_n})\in \mathbb{R}^{m_k\times m_1\times \cdots \times m_{k-1}\times m_{k+1}\times \cdots\times m_n}
$$
with
$$
\bar{a}^k_{i_1i_2\cdots i_n}=a^k_{i_ki_1\cdots i_{k-1}i_{k+1}\cdots i_n},\quad \forall i_j\in [m_j]\;\; \mbox{\rm and}\;\; j\in [n].
$$
\end{Definition}

Then, by Proposition \ref{prop-1}, we have $\bar{\mathscr{A}}^1=\mathscr{A}^1$; and for any
$x=\left((x^k)_{k\in [n]}\right)\in \mathbb{R}^m$, we have that
$$
\frac{\partial}{\partial x^k}\mathscr{A}^kx^1x^2\cdots x^n=\bar{\mathscr{A}}^kx^1\cdots x^{k-1}x^{k+1}\cdots x^n
$$
and
$$
\langle x^k,\bar{\mathscr{A}}^kx^1\cdots x^{k-1}x^{k+1}\cdots x^n\rangle=\mathscr{A}^kx^1x^2\cdots x^n
$$
hold for any $k\in [n]$.

Furthermore, by (\ref{E-new-tensor}), it is not difficult to see that
\begin{eqnarray}\label{E-add-1}
\mathscr{A}{x}^{m-1}=\left(\begin{array}{c}
                      \bar{\mathscr{A}}^1x^2\cdots x^n\\
                      \vdots \\
                      \bar{\mathscr{A}}^k x^1\cdots x^{k-1}x^{k+1}\cdots x^n\\
                      \vdots \\
                      \bar{\mathscr{A}}^n x^1x^2\cdots x^{n-1}
                      \end{array}\right).
\end{eqnarray}

Now, we can construct a tensor complementarity problem as follows:

Find $y=\left((y^k)_{k\in [n]}\right)\in \mathbb{R}^m$ such that
\begin{eqnarray}\label{app-tcp}
y\geq0,\quad \mathscr{A}y^{m-1}+q\geq0,\quad \langle y,\mathscr{A}y^{m-1}+q\rangle=0,
\end{eqnarray}
where $\mathscr{A}\in \mathbb{T}_{n,m}$ is a known tensor given by (\ref{E-new-tensor}), $q\in \mathbb{R}^m$ is a known vector given by
\begin{eqnarray*}
q:=\left(\begin{array}{c} -e_{m_1} \\ -e_{m_2} \\ \vdots \\ -e_{m_n} \end{array}\right)\in \mathbb{R}^{m_1}\times \mathbb{R}^{m_2}\times \cdots \times \mathbb{R}^{m_n}=\mathbb{R}^m,
\end{eqnarray*}
and $\mathscr{A}y^{m-1}$ is defined by (\ref{E-add-1}) by replacing $x$ by $y$.

\begin{Remark}
\begin{itemize}
\item[(i)] The constructed complementarity problem (\ref{app-tcp}) is a specific tensor complementarity problem. We denote the problem (\ref{app-tcp}) by the TCP$(q,\mathscr{A})$.
\item[(ii)] When $n=2$, tensors $\mathscr{A}^1$ and $\mathscr{A}^2$ reduce to two matrices, denoted by $A^1$ and $A^2$, respectively; and
the tensor $\mathscr{A}$ defined by (\ref{E-new-tensor}) reduces to a matrix $A$ given by
$$
A=\left(\begin{array}{cc} 0 & A^1 \\ {A^2}^T & 0 \end{array}\right).
$$
In this case, the TCP$(q,\mathscr{A})$ (\ref{app-tcp}) reduces to a linear complementarity problem, which is a reformulation of the bimatrix game \cite{lh-64}.
\end{itemize}
\end{Remark}

In the following, we will show that finding a Nash equilibrium point of the multilinear game is equivalent to finding a solution of the TCP$(q,\mathscr{A})$ (\ref{app-tcp}) with the explicit corresponding relation between the solutions of these two problems.

\begin{Theorem}\label{thm-main}
If $x^*=\left(({x^k}^*)_{k\in [n]}\right)\in \mathbb{R}^m$ is a Nash equilibrium point of the multilinear game, then $y^*=\left(({y^k}^*)_{k\in [n]}\right)\in \mathbb{R}^m$ defined by
\begin{eqnarray}\label{E-thm1-1}
{y^k}^*:=\sqrt[n-1]{\frac{(\mathscr{A}^k{x^1}^*{x^2}^*\cdots {x^n}^*)^{n-2}}{\prod_{i\in [n]\setminus \{k\}}\mathscr{A}^i{x^1}^*{x^2}^*\cdots {x^n}^*}}\;{x^k}^*\;\; \mbox{\rm for any}\;\; k\in [n]
\end{eqnarray}
is a solution of the TCP$(q,\mathscr{A})$ (\ref{app-tcp}).

Conversely, if $y^*=\left(({y^k}^*)_{k\in [n]}\right)\in \mathbb{R}^m$ is a solution of the TCP$(q,\mathscr{A})$ (\ref{app-tcp}), then ${y^k}^*\neq 0$ for any $k\in [n]$; and $x^*=\left(({x^k}^*)_{k\in [n]}\right)\in \mathbb{R}^m$ defined by
\begin{eqnarray}\label{E-thm1-2}
{x^k}^*:=\frac{{y^k}^*}{e^T_{m_k}{y^k}^*}\;\; \mbox{\rm for any}\;\;  k\in [n]
\end{eqnarray}
is a Nash equilibrium point of the multilinear game.
\end{Theorem}

\noindent {\bf Proof}. ``$\Longrightarrow$". Suppose that $x^*=\left(({x^k}^*)_{k\in [n]}\right)\in \mathbb{R}^m$ is a Nash equilibrium point of the multilinear game, we show that $y^*=\left(({y^k}^*)_{k\in [n]}\right)\in \mathbb{R}^m$ defined by (\ref{E-thm1-1}) is a solution of the TCP$(q,\mathscr{A})$ (\ref{app-tcp}). For any $k\in [n]$, by the KKT conditions of problem (\ref{E-opt-1}), there exist a number $\lambda_k^*\in \mathbb{R}$ and a nonnegative vector $\mu_k^*\in \mathbb{R}^{m_k}$ such that
\begin{eqnarray}\label{E-thm-1}
\bar{\mathscr{A}}^k{x^1}^*\cdots {x^{k-1}}^*{x^{k+1}}^*\cdots {x^n}^*-\lambda_k^* e_{m_k}-\mu_k^*=0
\end{eqnarray}
and
\begin{eqnarray}\label{E-thm-2}
e_{m_k}^T{x^k}^*=1, \quad  {x^k}^*\geq0, \quad \mu_k^*\geq 0,\quad {\mu_k^*}^T{x^k}^*=0.
\end{eqnarray}

By (\ref{E-thm-1}), it is easy to obtain that for any $k\in [n]$,
$$
\mathscr{A}^k{x^1}^*{x^2}^*\cdots {x^n}^*-\lambda_k^*e_{m_k}^T{x^k}^*-{\mu_k^*}^T{x^k}^*=0,
$$
which, together with equalities given in (\ref{E-thm-2}), implies that
\begin{eqnarray}\label{E-lambda-k}
\mathscr{A}^k{x^1}^*{x^2}^*\cdots {x^n}^*=\lambda_k^*e_{m_k}^T{x^k}^*-{\mu_k^*}^T{x^k}^*=\lambda_k^*.
\end{eqnarray}
Since ${x^k}^*\geq 0$ and ${x^k}^*\neq0$ for any $k\in [n]$; and $a^k_{i_1i_2\cdots i_n}>0$ for any $k\in [n]$ and any $i_j\in [m_j]$ with all $j\in [n]$, it is easy to show that
$$
\lambda_k^*=\mathscr{A}^k{x^1}^*{x^2}^*\cdots {x^n}^*>0,\quad \forall k\in [n].
$$
Thus, for any $k\in [n]$,
\begin{eqnarray}\label{non-neg}
{y^k}^*=\sqrt[n-1]{\frac{(\lambda_k^*)^{n-2}}{\prod_{i\in [n]\setminus \{k\}}\lambda_i^*}}\;{x^k}^*\geq 0.
\end{eqnarray}
Furthermore,
\begin{eqnarray}\label{thm1-add1}
\mathscr{A}{y^*}^{m-1}+q &=& \left(\begin{array}{c}
                      \bar{\mathscr{A}}^1{y^2}^*\cdots {y^n}^*-e_{m_1}\\
                      \vdots \\
                      \bar{\mathscr{A}}^k {y^1}^*\cdots {y^{k-1}}^*{y^{k+1}}^*\cdots {y^n}^*-e_{m_k}\\
                      \vdots \\
                      \bar{\mathscr{A}}^n {y^1}^*{y^2}^*\cdots {y^{n-1}}^*-e_{m_n}
                      \end{array}\right) \nonumber\\
                    &=&\left(\begin{array}{c}
                      \frac{1}{\lambda_1^*}\bar{\mathscr{A}}^1{x^2}^*\cdots {x^n}^*-e_{m_1}\\
                      \vdots \\
                      \frac{1}{\lambda_k^*}\bar{\mathscr{A}}^k {x^1}^*\cdots {x^{k-1}}^*{x^{k+1}}^*\cdots {x^n}^*-e_{m_k}\\
                      \vdots \\
                      \frac{1}{\lambda_n^*}\bar{\mathscr{A}}^n {x^1}^*{x^2}^*\cdots {x^{n-1}}^*-e_{m_n}
                      \end{array}\right) \nonumber\\
                   &=&\left(\begin{array}{c}
                      \frac{1}{\lambda_1^*}(\lambda_1^*e_{m_1}+\mu_1^*)-e_{m_1}\\
                      \vdots \\
                      \frac{1}{\lambda_k^*}(\lambda_k^*e_{m_k}+\mu_k^*)-e_{m_k}\\
                      \vdots \\
                      \frac{1}{\lambda_n^*}(\lambda_n^*e_{m_n}+\mu_n^*)-e_{m_n}
                      \end{array}\right)\nonumber\\
                   &=& \left(\begin{array}{c}
                      \frac{\mu_1^*}{\lambda_1^*}\\
                      \vdots \\
                      \frac{\mu_k^*}{\lambda_k^*}\\
                      \vdots \\
                      \frac{\mu_n^*}{\lambda_n^*}
                      \end{array}\right)\nonumber\\
                  &\geq& 0,
\end{eqnarray}
where the first equality follows from (\ref{E-add-1}), the second equality from (\ref{non-neg}), and the third equality and the last inequality from (\ref{E-thm-1}) and (\ref{E-thm-2}).
Moreover,
\begin{eqnarray}\label{thm1-add2}
{y^*}^T(\mathscr{A}{y^*}^{m-1}+q)
 &=& \left(\begin{array}{c}
                      {y^1}^*\\
                      \vdots \\
                      {y^k}^*\\
                      \vdots \\
                      {y^n}^*
                      \end{array}\right)^T
                     \left(\begin{array}{c}
                      \bar{\mathscr{A}}^1{y^2}^*\cdots {y^n}^*-e_{m_1}\\
                      \vdots \\
                      \bar{\mathscr{A}}^k {y^1}^*\cdots {y^{k-1}}^*{y^{k+1}}^*\cdots {y^n}^*-e_{m_k}\\
                      \vdots \\
                      \bar{\mathscr{A}}^n {y^1}^*{y^2}^*\cdots {y^{n-1}}^*-e_{m_n}
                      \end{array}\right) \nonumber\\
  &= & \sum\limits_{k=1}^n  {{y^k}^*}^T(\bar{\mathscr{A}}^k{y^1}^*\cdots {y^{k-1}}^*{y^{k+1}}^*\cdots {y^n}^*-e_{m_k}) \nonumber\\
  &=& \sum\limits_{k=1}^n \left\{\sqrt[n-1]{\frac{1}{\prod_{i=1}^n\lambda_i}}\mathscr{A}^k{x^1}^*\cdots {x^n}^*-\sqrt[n-1]{\frac{(\lambda_k^*)^{n-2}}{\prod_{i\in [n]\setminus\{k\}}\lambda_i^*}}e_{m_k}^T{x^k}^*\right\} \nonumber\\
  &= & \sum\limits_{k=1}^n \sqrt[n-1]{\frac{(\lambda_k^*)^{n-2}}{\prod_{i\in [n]\setminus\{k\}}\lambda_i^*}}\;(1-e_{m_k}^T{x^k}^*) \nonumber\\
  &= & 0,
\end{eqnarray}
where the third equality holds by (\ref{non-neg}), the forth equality holds by (\ref{E-lambda-k}), and the last equality holds by (\ref{E-thm-2}).

Combining (\ref{non-neg}) with (\ref{thm1-add1}) and (\ref{thm1-add2}), we obtain that $y^*=\left(({y^k}^*)_{k\in [n]}\right)\in \mathbb{R}^m$ defined by (\ref{E-thm1-1}) is a solution of the TCP$(q,\mathscr{A})$ (\ref{app-tcp}).

``$\Longleftarrow$". Suppose that $y^*=\left(({y^k}^*)_{k\in [n]}\right)\in \mathbb{R}^m$ is a solution of the TCP$(q,\mathscr{A})$ (\ref{app-tcp}), then
\begin{eqnarray}\label{formula}
\begin{array}{l}
\left(\begin{array}{c}
                      {y^1}^*\\
                      \vdots \\
                      {y^k}^*\\
                      \vdots \\
                      {y^n}^*
                      \end{array}\right)\geq 0, \quad
\left(\begin{array}{c}
                      \bar{\mathscr{A}}^1{y^2}^*\cdots {y^n}^*-e_{m_1}\\
                      \vdots \\
                      \bar{\mathscr{A}}^k {y^1}^*\cdots {y^{k-1}}^*{y^{k+1}}^*\cdots {y^n}^*-e_{m_k}\\
                      \vdots \\
                      \bar{\mathscr{A}}^n {y^1}^*{y^2}^*\cdots {y^{n-1}}^*-e_{m_n}
                      \end{array}\right)\geq 0, \vspace{2mm}\\
\left(\begin{array}{c}
                      {y^1}^*\\
                      \vdots \\
                      {y^3}^*\\
                      \vdots \\
                      {y^n}^*
                      \end{array}\right)^T
                     \left(\begin{array}{c}
                      \bar{\mathscr{A}}^1{y^2}^*\cdots {y^n}^*-e_{m_1}\\
                      \vdots \\
                      \bar{\mathscr{A}}^k {y^1}^*\cdots {y^{k-1}}^*{y^{k+1}}^*\cdots {y^n}^*-e_{m_k}\\
                      \vdots \\
                      \bar{\mathscr{A}}^n {y^1}^*{y^2}^*\cdots {y^{n-1}}^*-e_{m_n}
                      \end{array}\right)=0.
\end{array}
\end{eqnarray}
It is easy to show that ${y^k}^*\neq 0$ for any $k\in [n]$. In fact, if ${y^k}^*=0$ for some $k\in [n]$, then by the second inequality in (\ref{formula}), we have that $-e_{m_j}\geq 0$ for any $j\in [n]\setminus \{k\}$, which is a contradiction.

Next, we prove that $x^*=\left(({x^k}^*)_{k\in [n]}\right)\in \mathbb{R}^m$ defined by (\ref{E-thm1-2}) is a Nash equilibrium point of the multilinear game. For this purpose, we need to show that, for any $k\in [n]$, there exist a number $\lambda_k^*\in \mathbb{R}$ and a nonnegative vector $\mu_k^*\in \mathbb{R}^{m_k}$ such that
\begin{eqnarray}\label{E-thm-1-1}
\bar{\mathscr{A}}^k{x^1}^*\cdots {x^{k-1}}^*{x^{k+1}}^*\cdots {x^n}^*-\lambda_k^* e_{m_k}-\mu_k^*=0
\end{eqnarray}
and
\begin{eqnarray}\label{E-thm-2-1}
e_{m_k}^T{x^k}^*=1, \quad  {x^k}^*\geq0, \quad \mu_k^*\geq 0,\quad {\mu_k^*}^T{x^k}^*=0.
\end{eqnarray}

By (\ref{formula}), we have that for any $k\in [n]$,
$$
{{y^k}^*}^T(\bar{\mathscr{A}}^k{y^1}^*\cdots {y^{k-1}}^*{y^{k+1}}^*\cdots {y^n}^*-e_{m_k})=0,
$$
i.e.,
$$
\mathscr{A}^k{y^1}^*{y^2}^*\cdots {y^n}^*-e_{m_k}^T{y^k}^*=0.
$$
For any $k\in [n]$, since ${y^k}^*\neq 0$ and ${y^k}^*\geq 0$, we have that  $e_{m_k}^T{y^k}^*>0$; and then
\begin{eqnarray*}
\mathscr{A}^k\frac{{y^1}^*}{e_{m_1}^T{y^1}^*}\frac{{y^2}^*}{e_{m_2}^T{y^2}^*}\cdots \frac{{y^n}^*}{e_{m_n}^T{y^n}^*}-\frac{1}{\prod_{i\in [n]\setminus \{k\}} e_{m_i}^T{y^i}^*}=0.
\end{eqnarray*}
By (\ref{E-thm1-2}), the above equality becomes
\begin{eqnarray}\label{formula3}
\mathscr{A}^k{x^1}^*{x^2}^*\cdots {x^n}^*-\frac{1}{\prod_{i\in [n]\setminus \{k\}} e_{m_i}^T{y^i}^*}=0.
\end{eqnarray}
For any $k\in [n]$, from ${y^k}^*\geq 0$, $e_{m_k}^T{y^k}^*>0$ and the definition of ${x^k}^*$, it follows that ${x^k}^*\geq 0$ and $e_{m_k}^T{x^k}^*=1$.
In addition, for any $k\in [n]$, since
$$
\bar{\mathscr{A}}^k{y^1}^*\cdots {y^{k-1}}^*{y^{k+1}}^*\cdots {y^n}^*-e_{m_k}\geq 0,
$$
we have that
$$
\bar{\mathscr{A}}^k{x^1}^*\cdots {x^{k-1}}^*{x^{k+1}}^*\cdots {x^n}^*-\frac{e_{m_k}}{\prod_{i\in [n]\setminus \{k\}}e_{m_i}^T{y^i}^*}\geq 0,
$$
which implies that there exists a nonnegative vector $\mu_k^*\in \mathbb{R}^{m_k}$ such that
$$
\bar{\mathscr{A}}^k{x^1}^*\cdots {x^{k-1}}^*{x^{k+1}}^*\cdots {x^n}^*-\frac{e_{m_k}}{\prod_{i\in [n]\setminus \{k\}}e_{m_i}^T{y^i}^*}-\mu_k^* = 0;
$$
and furthermore,
\begin{eqnarray*}
{\mu_k^*}^T{x^k}^* &=& {{x^k}^*}^T\left(\bar{\mathscr{A}}^k{x^1}^*\cdots {x^{k-1}}^*{x^{k+1}}^*\cdots {x^n}^*-\frac{e_{m_k}}{\prod_{i\in [n]\setminus \{k\}}e_{m_i}^T{y^i}^*}\right) \\
  &=&  \mathscr{A}^k{x^1}^*{x^2}^*\cdots {x^n}^*-\frac{e_{m_k}^T{x^k}^*}{\prod_{i\in [n]\setminus \{k\}}e_{m_i}^T{y^i}^*} \\
  &=&  \mathscr{A}^k{x^1}^*{x^2}^*\cdots {x^n}^*-\frac{1}{\prod_{i\in [n]\setminus \{k\}}e_{m_i}^T{y^i}^*} \\
  &=&  0,
\end{eqnarray*}
where the last equality holds by (\ref{formula3}). So, we obtain that (\ref{E-thm-1-1}) and (\ref{E-thm-2-1}) holds with
$$
\lambda_k^*=\frac{1}{\prod_{i\in [n]\setminus \{k\}}e_{m_i}^T{y^i}^*}.
$$
Therefore, $x^*=\left(({x^k}^*)_{k\in [n]}\right)\in \mathbb{R}^m$ defined by (\ref{E-thm1-2}) is a Nash equilibrium point of the multilinear game.
\ep
\begin{Remark}
\begin{itemize}
\item[(i)] In Theorem \ref{thm-main}, we have reformulated the multilinear game as a tensor complementarity problem; and especially, we have established a one to one correspondence between the solutions of these two classes of problems, which built a bridge between two classes of problems.
\item[(ii)] When $n=2$, the TCP$(q,\mathscr{A})$ (\ref{app-tcp}) reduces to a linear complementarity problem, which is a reformulation of the bimatrix game \cite{hxq-06,lh-64,Isac-00}; and the results of Theorem \ref{thm-main} reduces to those obtained in the case of the bimatrix game \cite{hxq-06,Isac-00}.
\item[(iii)] By using Theorem \ref{thm-main}, we can investigate the TCP$(q,\mathscr{A})$ (\ref{app-tcp}) by the known results on the $n$-person noncooperative game. It is easy to see that the multilinear game has at least a Nash equilibrium point by using Nash's result \cite{nash-51}; and hence, by Theorem \ref{thm-main}, we obtain that the TCP$(q,\mathscr{A})$ (\ref{app-tcp}) has at least a solution.
\item[(iv)] By using Theorem \ref{thm-main}, we can also investigate the multilinear game by using the theory and methods for the nonlinear complementarity problems \cite{fp-03,hxq-06}. In the next section, we apply a smoothing-type algorithm to solve the TCP$(q,\mathscr{A})$ (\ref{app-tcp}).
\end{itemize}
\end{Remark}

\section{Algorithm and numerical results}
\hspace{4mm} It is well known that the smoothing-type algorithm is a class of effective methods for solving variational inequalities and  complementarity problems, and related optimization problems \cite{bx-00,ch-93,cqs-98,fl-04,flt-01,hn-10,qsz-00}. In this section, we apply a smoothing-type algorithm to solve the TCP$(q,\mathscr{A})$ (\ref{app-tcp}) and give some preliminary numerical results for solving the multilinear games.

Let payoff tensors of the multilinear game be given by $\mathscr{A}^k$ for all $k\in [n]$, the tensor $\mathscr{A}\in \mathbb{T}_{n,m}$ be defined by (\ref{E-new-tensor}), and the tensors $\bar{\mathscr{A}}^k$ for all $k\in [n]$ be defined by Definition \ref{def-tensor-bar-a}. Denote
$$
F(y):=\left(\begin{array}{c}
                      \bar{\mathscr{A}}^1{y^2}^*\cdots {y^n}^*-e_{m_1}\\
                      \vdots \\
                      \bar{\mathscr{A}}^k {y^1}^*\cdots {y^{k-1}}^*{y^{k+1}}^*\cdots {y^n}^*-e_{m_k}\\
                      \vdots \\
                      \bar{\mathscr{A}}^n {y^1}^*{y^2}^*\cdots {y^{n-1}}^*-e_{m_n}
                      \end{array}\right),
$$
then we can rewrite the TCP$(q,\mathscr{A})$ (\ref{app-tcp}) as follows: Find $y=\left((y^k)_{k\in [n]}\right)\in \mathbb{R}^m$ and $s=\left((s^k)_{k\in [n]}\right)\in \mathbb{R}^m$ such that
\begin{eqnarray}\label{app-tcp-1}
y\geq0,\quad s=F(y)\geq0,\quad \langle y,s\rangle=0.
\end{eqnarray}
We define a function $H: \mathbb{R}^{1+2m}\rightarrow \mathbb{R}^{1+2m}$ by
\begin{eqnarray*}
H(\mu,y,s):=\left(\begin{array}{c} \mu \\ s-F(y) \\ \Phi(\mu,y,s)+\mu y \end{array}\right),
\end{eqnarray*}
where $\Phi(\mu,y,s)=(\phi(\mu,y_1,s_1),\phi(\mu,y_2,s_2),\ldots,\phi(\mu,y_m,s_m))^T$ with
$$
\phi(\mu,y_i,s_i)=y_i+s_i-\sqrt{(y_i-s_i)^2+4\mu}, \forall i\in \{1,2,\ldots,m\}.
$$
It is obvious that $(y,s)$ solves the problem (\ref{app-tcp-1}) if and only if $H(\mu,y,s)=0$. Since the function $H$ is continuously differentiable for any $(\mu,y,s)\in \mathbb{R}^{1+2m}$ with $\mu>0$, we can apply some Newton-type methods to solve the system of smooth equations $H(\mu,y,s)=0$ at each iteration and make $\mu\rightarrow 0$ so that a solution of the problem (\ref{app-tcp-1}) can be found. We use the following algorithm to solve the problem (\ref{app-tcp-1}).
\begin{Algorithm}\label{algo}(A Smoothing-type Algorithm)
\begin{description}
\item [Step 0] Choose $\delta, \sigma\in (0,1)$. Let $\mu_0>0$ and $(y^0,s^0)\in \mathbb{R}^{2m}$ be an arbitrary vector. Set $z^0:=(\mu_0,y^0,s^0)$. Choose $\beta>1$ such that $\|H(z^0)\|\le \beta\mu_0$. Set $e^0:=(1,0,\ldots,0)\in \mathbb{R}^{1+2m}$ and $k:=0$.

\item [Step 1]  If $\|H(z^k)\|=0$, stop.

\item [Step 2]  Compute $\d z^k:=(\d \mu_k,\d x^k,\d s^k)\in \mathbb{R}\times \mathbb{R}^m \times \mathbb{R}^{m}$ by
\begin{eqnarray}\label{newton-equa}
H(z^k)+ H^\prime (z^k)\d z^k=(1/\beta) \|H(z^k)\|e^0.
\end{eqnarray}

\item [Step 3] Let $\lambda_k$ be the maximum of the values $1,\delta,\delta^2,\cdots$ such that
\begin{eqnarray}\label{linesearch}
\|H(z^k+\lambda_k\d z^k)\|\le [1-\sigma (1-1/\beta)\lambda_k]\|H(z^k)\|.
\end{eqnarray}

\item [Step 4]  Set $z^{k+1}:=z^k+\lambda_k\d z^k$ and $k:=k+1$. Go to Step 1.
\end{description}
\end{Algorithm}

The above algorithmic framework was proposed in \cite{huang-05}; and from \cite{huang-05} it follows that Algorithm \ref{algo} is globally convergent under suitable assumptions.

In the following, we give some preliminary numerical results of Algorithm \ref{algo} for solving the multilinear game and bimatrix game. Throughout our experiments, the parameters used in Algorithm \ref{algo} are chosen as
$$
\delta:=0.75,\quad \sigma:=10^{-4},\quad y^0:=0.01*ones(m,1),\quad s^0:=F(y^0),
$$
where $m$ is given in the tested examples. In our experiments, we take $\mu_0:=0.1$; and if the algorithm fails to find a solution to the TCP$(q,\mathscr{A})$ (\ref{app-tcp}), we try to take $\mu_0:=0.01$ or $\mu_0:=0.1+3*p$ for $p=2,3,4,5,6$, respectively.
We denote $z^0:=(\mu_0,y^0,s^0)$ and take $\beta := \|H(z^0)\|/\mu_0$. We use $\|H(z^k)\|\le 10^{-6}$ as the stopping rule.

\begin{Example}\label{exam0}
Consider the multilinear games with three players, where the payoff tensors $\mathscr{A}^1, \mathscr{A}^2,\mathscr{A}^3\in \mathbb{R}^{m_1\times m_2\times m_3}$ are randomly generated by using $rand(m_1,m_2,m_3)$, respectively.
\end{Example}

Obviously, for different values of $m_1$, $m_2$ and $m_3$, different games are generated by Example \ref{exam0}. We use Algorithm \ref{algo} to solve these games, where values of $m_1$, $m_2$, $m_3$ and $m:=m_1+m_2+m_3$ are specified in the table of the numerical results. In our experiments, for any fixed $m_1$, $m_2$ and $m_3$, the random problems are generated ten times for which Algorithm \ref{algo} can find an approximation solution to every generated problem. The numerical results are listed in Table 1, where AI (MinI and MaxI) denotes the average number (minimal number and maximal number) of iterations for solving ten randomly generated problems of each size; AT (MinT and MaxT) denotes the average (minimal and maximal) CPU time in second for solving ten randomly generated problems of each size; and ARes denotes the average value of $\|H(z^k)\|$ for ten randomly generated problems of each size when the algorithm stops.

\begin{table}[ht] \label{table1}
  \caption{The numerical results of the problem in Example \ref{exam0}}
  \begin{center}
    \begin{tabular}[c]
      {| c | c | c | c | c | c | c|}
      \hline
      $m$  & $m_1$ & $m_2$ & $m_3$  & AI/MinI/MaxI    & AT/MinT/MaxT(s) & ARes\\
      \hline
      \hline
           & 2     & 2     & 6      & 13.7/9/23       &0.0546/0.0156/0.125    & $2.68e-7$   \\
           & 2     & 3     & 5      & 16.3/12/31      &0.0764/0.0156/0.187    & $1.17e-7$   \\
           & 2     & 4     & 4      & 15.5/10/26      &0.0889/0.0156/0.265    & $1.55e-7$   \\
           & 3     & 5     & 2      & 15.1/10/22      &0.125/0.0468/0.281     & $1.98e-7$   \\
        10 & 3     & 2     & 5      & 13.9/9/24       &0.0842/0.0156/0.203    & $3.35e-7$   \\
           & 4     & 4     & 2      & 18.3/11/25      &0.0889/0.0156/0.265    & $1.45e-7$   \\
           & 4     & 2     & 4      & 13.6/9/22       &0.0827/0.0156/0.172    & $3.35e-7$   \\
           & 5     & 3     & 2      & 15.0/9/21       &0.0515/0.0156/0.140    & $2.64e-7$   \\
           & 6     & 2     & 2      & 12.7/9/26       &0.0250/0.0156/0.0468   & $2.12e-7$   \\ \hline
           & 3     & 5     & 12     & 23.8/17/33      &0.200/0.0624/0.421    & $2.04e-7$   \\
           & 3     & 8     & 9      & 21.8/14/31      &0.133/0.0312/0.328    & $2.51e-7$   \\
           & 3     & 12    & 5      & 22.5/13/29      &0.179/0.0936/0.472    & $2.34e-7$   \\
           & 4     & 6     & 10     & 24.0/16/35      &0.137/0.0312/0.281    & $2.44e-7$   \\
        20 & 4     & 8     & 8      & 22.5/14/32      &0.136/0.0468/0.437    & $3.08e-7$   \\
           & 4     & 10    & 6      & 21.6/14/39      &0.179/0.0312/0.374    & $2.38e-7$   \\
           & 6     & 5     & 9      & 21.0/13/40      &0.125/0.0156/0.328    & $1.30e-7$   \\
           & 8     & 4     & 8      & 27.2/13/39      &0.183/0.0312/0.484    & $3.29e-7$   \\
           & 12    & 5     & 3      & 22.1/14/32      &0.212/0.0780/0.593    & $3.04e-7$   \\ \hline
    \end{tabular}
  \end{center}
\end{table}

From Table 1, it is easy to see that the TCP$(q,\mathscr{A})$ (\ref{app-tcp-1}) can be effectively solved by
Algorithm \ref{algo}. Furthermore, by Theorem \ref{thm-main} we can obtain that a Nash equilibrium point of the concerned game can be found by using Algorithm \ref{algo}. In order to see this more clearly, we test two specific problems in the following.

\begin{Example}\label{exam1}
Consider a multilinear game with three players, where the payoff tensors $\mathscr{A}^1,\mathscr{A}^2,\mathscr{A}^3\in \mathbb{R}^{2\times 3\times 2}$ are given by
\begin{eqnarray*}
\begin{array}{l}
\mathscr{A}^1(:,:,1)=\left(\begin{array}{ccc}
    0.0605 &   0.5269  &  0.6569 \\
    0.3993 &   0.4168  &  0.6280
    \end{array}\right), \quad
    \mathscr{A}(:,:,2)=\left(\begin{array}{ccc}
    0.2920  &  0.0155  &  0.1672 \\
    0.4317  &  0.9841  &  0.1062
    \end{array}\right), \\
\mathscr{A}^2(:,:,1)=\left(\begin{array}{ccc}
    0.3724  &  0.4897  &  0.9516 \\
    0.1981  &  0.3395  &  0.9203
    \end{array}\right), \quad
    \mathscr{A}^2(:,:,2)=\left(\begin{array}{ccc}
    0.0527  &  0.2691 &   0.5479 \\
    0.7379  &  0.4228 &   0.9427
    \end{array}\right), \\
\mathscr{A}^3(:,:,1)=\left(\begin{array}{ccc}
    0.4177  &  0.3015 &   0.6663 \\
    0.9831  &  0.7011 &   0.5391
    \end{array}\right), \quad
    \mathscr{A}^2(:,:,2)=\left(\begin{array}{ccc}
    0.6981  &  0.1781  &  0.9991 \\
    0.6665  &  0.1280  &  0.1711
    \end{array}\right).
\end{array}
\end{eqnarray*}
\end{Example}

We use Algorithm \ref{algo} to solve the TCP$(q,\mathscr{A})$ (\ref{app-tcp-1}) with the payoff tensors being given by Example \ref{exam1}; and a solution to this tensor complementarity problem:
\begin{eqnarray*}
\begin{array}{l}
y^*=(0.6235,0.0000,3.8396,0.0000,0.0000,4.3070,0.0000)^T,\\
s^*=(0.0000,5.6024,0.0000,0.3149,1.5553,0.0000,0.6711)^T
\end{array}
\end{eqnarray*}
is obtained with $10$ iterative steps in $0.0156$ seconds. Furthermore, by Theorem \ref{thm-main} we obtain that a Nash equilibrium point of the concerned game is $x^*=({x^1}^*,{x^2}^*,{x^3}^*)$ with
$$
{x^1}^*=(1.0000,0.0000)^T,\quad {x^2}^*=(1.0000,0.0000,0.0000)^T,\quad {x^3}^*=(1.0000,0.0000)^T.
$$

\begin{Example}\label{exam2}
Consider the bimatrix game ``Battle of the Sexes" \cite{ms-64}, where two payoff matrices ${A}^1,{A}^2\in \mathbb{R}^{2\times 2}$ are given by
\begin{eqnarray*}
{A}^1=\left(\begin{array}{cc}  2 & -1 \\ -1 & 1
    \end{array}\right), \quad
{A}^2=\left(\begin{array}{ccc} 1 & -1 \\ -1 & 2
    \end{array}\right).
\end{eqnarray*}
\end{Example}

We use Algorithm \ref{algo} to solve the TCP$(q,\mathscr{A})$ (\ref{app-tcp-1}) with the payoff matrices being given by Example \ref{exam2}; and a solution to this tensor complementarity problem:
\begin{eqnarray*}
y^*=(3,2,2,3)^T,\quad s^*=(0,0,0,0)^T
\end{eqnarray*}
is obtained with $5$ iterative steps in $0.0156$ seconds. Furthermore, by Theorem \ref{thm-main} we obtain that a Nash equilibrium point of the concerned game is $x^*=({x^1}^*,{x^2}^*)$ with
$$
{x^1}^*=(0.6,0.4)^T,\quad {x^2}^*=(0.4,0.6)^T.
$$

From these numerical results, we can see that Algorithm \ref{algo} is effective for solving the tensor complementarity problem (\ref{app-tcp-1}). We have also tested some other problems, the computation effect is similar.

\section{Conclusions}
\hspace{4mm}
In this paper, we reformulated the multilinear game as a tensor complementarity problem and showed that finding a Nash equilibrium point of the multilinear game is equivalent to finding a solution of the resulted tensor complementarity problem. Especially, we provided a one to one correspondence between the solutions of the multilinear game and the tensor complementarity problem, which built a bridge between these two classes of problems so that one can investigate one problem by using the theory and methods for another problem. We also applied a smoothing-type algorithm to solve the resulted tensor complementarity problem and reported some preliminary numerical results for solving the multilinear games. Hopefully some more effective algorithms can be designed to solve the tensor complementarity problem by using the structure of the tensors and properties of the homogeneous polynomials.

\end{document}